\documentclass[a4paper,landscape]{article}

\usepackage{a4} 
\usepackage{amsfonts}
\usepackage{amssymb}
\usepackage{amsmath}
\usepackage{amscd}
\usepackage{amsthm}
\usepackage{epsfig}
\usepackage[all]{xy}
\usepackage{color}
\usepackage{bbm}
\usepackage{graphicx}
\usepackage{psfrag}

\newcommand{\R}{\mathbb{R}}

\newcommand{\Z}{\mathbb{Z}}

\newcommand{\Le}{{\cal L}}
\newcommand{\Ka}{{\sf K}}
\newcommand{\La}{{\sf L}}

\newcommand{\ra}{\rightarrow}

\newcommand{\bprf}{\begin{proof}[Proof]}
\newcommand{\eprf}{\end{proof}}

\newtheorem{thm}{Theorem}

\newtheorem{conj}[thm]{Conjecture}

\newtheorem{prop}[thm]{Proposition}
\newtheorem{lem}[thm]{Lemma}

\begin{document}
\title{On the number of Tverberg \\ partitions in the prime power case}
\author{{\Large Stephan Hell}
\thanks{This research was supported by the
    Deutsche Forschungsgemeinschaft within the European graduate
    program `Combinatorics, Geometry, and Computation' (No. GRK
    588/2).}} 
\date{Institut f\"ur Mathematik, MA 6--2,
TU Berlin,\\ D--10623 Berlin, Germany,
hell@math.tu-berlin.de}

\maketitle

\begin{abstract}
We give an extension of the lower bound of \cite{vz93} for the
number of Tverberg partitions from the prime to the prime power
case. Our proof is inspired by the $\Z_p$--index version of the proof 
in \cite{mat03} and uses Volovikov's Lemma. Analogously, one obtains
an extension of the lower bound for the number of different
splittings of a generic necklace to the prime power case.
\end{abstract}
\begin{section}{Introduction}\label{sec-intro}
In 1966, Helge Tverberg showed that any set of $(d+1)(q-1)+1$ points
in $\R^d$ admits a partition into $q$ subsets such that the
intersection of their convex hulls is non--empty. Such partitions are
called Tverberg partitions; the result is best possible: For less than
$(d+1)(q-1)+1$ points in $\R^d$ the implication of the statement does 
not hold. Moreover,
it can be formulated in the following way.
\begin{thm}[\cite{tve66}]\label{thm-tverberg}Let $q\geq 2$, $d\geq 1$, 
and put $N:=(d+1)(q-1)$. For every affine map $f:\|\sigma^N\|\rightarrow\mathbb{R}^d$ 
there are $q$ disjoint faces $F_1,F_2,\ldots,F_q$ of the 
standard $N$--simplex $\sigma^N$ whose images under $f$ intersect: 
$\bigcap_{i=1}^{q}f(\|F_i\|)\not=\emptyset$.
\end{thm}
Relaxing affine maps to continuous maps one gets a more general 
problem which is known as the Topological Tverberg Theorem. For $q$
a prime this topological version was first proved by B\'ar\'any et
al.~\cite{bss81}. The proof uses a Borsuk--Ulam type argument and can
be found in  Matou\v{s}ek's book \cite{mat03} on topological methods
in combinatorics and geometry.
In 1987, \"Ozaydin proved the case $q$ being a prime power in an
unpublished manuscript \cite{oez87}, later Volovikov gave another
proof in \cite{vol96}. Both proofs make use of deep results from
algebraic topology. For arbitrary $q$ the problem is still open.

Theorem~\ref{thm-tverberg} establishes the existence of Tverberg
partitions. Another natural question is to ask for a lower bound:~How
many Tverberg partitions into $q$ subsets are there for a chosen
affine or continuous map $f$? Sierksma conjectured that there are at
least $((q-1)!)^d$ for any set of $(d+1)(q-1)+1$ points in $\R^d$. The
conjecture is still not proved. The case $d=1$ and
arbitrary $q$ can be proved for continuous maps using the intermediate
value theorem. The only non--trivial lower bound is established for $q$
being prime using a Borsuk--Ulam type argument (see \cite{vz93}). 
The following extends the result of \cite{vz93} to the prime 
power case using Volovikov's lemma from \cite{vol96}.
\begin{thm}\label{thm-main}Let $q=p^r$ be a prime power. For any
  continuous map \mbox{$f:\|\sigma^N\|\ra \R^d$}, where
  $N=(d+1)(q-1)$, the number of unordered $q$--tuples
  $\{F_1,F_2,\ldots,F_q\}$ of disjoint faces of the $N$--simplex with
  $\bigcap_{i=1}^{q}f(\| F_i \|)\not=\emptyset$ is at least
\[\frac{1}{(q-1)!}\cdot\left(\frac{q}{r+1}\right)^{\lceil\frac{N}{2}\rceil}\]
\end{thm}
A simplified proof for the lower bound of \cite{vz93} can be found 
in Section 6.6 of \cite{mat03}. In the prime power case $q=p^r$, we cannot use 
the $\Z_q$--action by cyclic shifting of the $q$ coordinates of the $q$--fold 
join as the space 
$(\R^{d})^{*q}_{\Delta}$ is a non--free $\Z_q$--space so that
$\mbox{ind}_{\Z_q}((\R^{d})^{*q}_{\Delta})=+\infty$.

\begin{table*}[h]
\begin{center}
\begin{tabular}{l|l|l|l}
Lower bound $\setminus$ $q$ & prime & prime power & arbitrary \\ \hline
1 & \cite{bss81} & \cite{oez87},\cite{vol96} & open\\ \hline
\cite{vz93}--type & \cite{vz93} & \checkmark& open \\ \hline
Sierksma & open & open & open 
\end{tabular}
\caption{Current state around the Topological Tverberg Theorem}
\end{center}
\end{table*}

Progress towards the general case has been slow. But recently 
T.~Sch\"oneborn \cite{sch04} was able to connect the 
Topological Tverberg Theorem to geometric graph theory type 
questions. In particular, he showed that the $d=2$ case is
equivalent to the following conjecture.
\begin{conj}[Winding partitions]\label{winding}
For every drawing of the complete graph $K_{3(q-1)+1}$ there
are either $q-1$ disjoint triangles of edges and a vertex $v$
or $q-2$ disjoint triangles of edges and an intersection point 
$p$ of two edges such that the winding number about $v$ 
resp.~about $p$ of each triangle is non--zero. 
\end{conj}
Here a drawing of a graph $G$ is a continuous map from $G$, seen as a 
one--dimensional simplicial complex, to the plane such that (i) no two 
vertices coincide, (ii) no edge passes through a point (except its
endpoints), (iii) no three edges intersect (outside their endpoints).
Any lower bound for the number of Tverberg partitions carries over to
winding partitions.

We give a proof of Theorem~\ref{thm-main} in Section~\ref{sec-tow}. In 
Section~\ref{sec-necklace} we sketch how to extend the lower bound for  
splitting generic necklaces of \cite{vz93} to the prime power case. 
I thank Juliette Blanca, Mark de Longueville, Ji\v{r}\'{i}
Matou\v{s}ek, Torsten Sch\"oneborn and G\"unter M.~Ziegler for 
helpful discussions and remarks.
\end{section}
\begin{section}{Preliminaries}\label{sec-prel}
Before proving our lower bound we repeat some definitions and
results from \cite{mat03}, mainly for fixing our notation. 
We write $[n]$ for the set $\{1,2,\ldots,n\}$.
Let $G$ be a finite group. A topological
space $X$ equipped with a (left) $G$--action $\Phi:G\ra \mbox{Homeo}(X)$ is
called a {\it $G$--space}; we write $g\, x$ for
$\Phi(g)(x)$. Continuous maps between $G$--spaces $X$ and $Y$
that respect the $G$--actions of $X$ and $Y$ are called {\it $G$--maps} or
{\it equivariant maps}. For $x\in X$ the set $O_x=\{g\, x\,|\,g\in G\}$ is
called the {\it orbit} of $x$. A $G$--space $(X,\Phi)$ where every $O_x$
has at least two elements is called {\it fixed point free}, i.~e.~no
point of X is fixed by all group elements. Let $X$ be a fixed point
free $G$--space and $Y\subset X$ closed under the $G$--action, then
$Y$ with the induced action of $X$ is again a fixed point free $G$--space.

The {\it join} $X*Y$ of spaces $X$ and $Y$ is a standard construction in
topology. One way of looking at it is to identify it with the set of
formal convex combinations \mbox{$tx\oplus (1-t)y$}, where $t\in [0,1],\,
x\in X,\, y\in Y$. We use the symbol $\oplus$ to underline that
the sum is formal and does not commute for $X=Y$. With this
identification the $n$--fold join $X^{*n}$ becomes the set of all
formal convex combinations \mbox{$t_1x_1\oplus t_2x_2\oplus\cdots\oplus
t_nx_n$}, where $t_1,t_2,\ldots,t_n$ are non--negative reals summing up
to $1$ and $x_1,x_2,\ldots,x_n\in X$. The join of simplicial complexes
is again a simplicial complex. For abstract simplicial complexes $\Ka$
and $\La$ the join is defined as the set of simplices $\{F\uplus
G\,|\,F\in\Ka,\, G\in \La\}$, where $F\uplus G=(F\times\{1\})
\cup (G\times\{2\})$ is the disjoint union of $F$ and $G$. For subsets
$A\subset\R^n$ and $B\subset\R^m$ of Euclidean spaces the join can be
represented {\it geometrically} in the following way: Embed $A\subset\R^n
\subset\R^{n+m+1}$ in the standard way, and embed
$B\subset\R^m\subset\R^{n+m+1}$ such that the first $n$ coordinates are equal
to $0$ and the last one is equal to $1$. The subspace
$C\subset\R^{n+m+1}$ defined as the union
of all segments joining a point of $A$ with a point of $B$ is
homeomorphic to $A*B$.
Finally, there is an inequality for the connectivity of the join $X*Y$ 
for topological spaces $X$ and $Y$: 
\begin{eqnarray}\label{conn}
\mbox{conn}(X*Y)\geq \mbox{conn}(X)+\mbox{conn}(Y)+2,
\end{eqnarray}
where a disconnected space has connectivity $-1$.

Let $n\geq k\geq 2$. We call an $n$--tuple $(x_1,x_2,\ldots,x_n)$ {\it
  $k$--wise distinct} if no $k$ among the $x_i$ are equal. 
The {\it $n$--fold $k$--wise deleted join} of a space $X$ is 
\[X^{*n}_{\Delta(k)}:=X^{*n}\setminus\{{\textstyle\frac{1}{n}x_1\oplus\frac{1}{n}x_2
\oplus \cdots\oplus\frac{1}{n}x_n\,}|\,(x_1,x_2,\ldots,x_n)
\mbox{ not $k$--wise distinct}\}.\]
In the case $k=n$ we delete the diagonal of $X^{*n}$, and 
for $k_1<k_2$ we have $X^{*n}_{\Delta(k_1)}\subset
X^{*n}_{\Delta(k_2)}$; we write $X^{*n}_{\Delta}$ for $X^{*n}_{\Delta(n)}$. 
For a simplicial complex $\Ka$ we define its {\it $n$--fold
  $k$--wise deleted join} as the following set of simplices:
 \[\Ka^{*n}_{\Delta(k)}:=\{ F_1\uplus F_2\uplus\cdots\uplus F_n\in
\Ka^{*n}\,|\, (F_1,F_2,\ldots,F_n) \mbox{ $k$--wise disjoint}\},\]
where an $n$--tuple $(F_1,F_2,\ldots,F_n)$ is called $k$--wise
disjoint if no $k$ among them have a non--empty intersection.
For simplicial complexes $\Ka$ we have 
$\|\Ka^{*n}_{\Delta(k)}\|\subset \|\Ka\|^{*n}_{\Delta(k)}$.
In the proof, we are interested in the special cases $k=2$ and $k=n$.

{\bf The group action}. The symmetric group $S_q$ acts on 
a (deleted) $q$--fold join by permuting the
$q$ coordinates. The following result is the key lemma in \cite{vol96}
for the prime power case $q=p^r$, and it is proved for actions of the subgroup
$G:=(\Z_p)^r$ of $S_q$. $G$ is a subgroup of $S_q$ in the
following way: number its $q$ elements in lexicographic order; 
an element $g\in G$ defines an isomorphism on $G$ by translation
$h\mapsto g+h$ for $h\in G$. Now the element $(1,1)\in (\Z_3)^2$ 
acts on $X^{*9}$:
\mbox{$t_1x_1\oplus\cdots\oplus t_9x_9\mapsto t_9x_9\oplus
  t_7x_7\oplus t_8x_8\oplus t_3x_3\oplus t_1x_1\oplus t_2x_2
\oplus t_6x_6\oplus t_4x_4\oplus t_5x_5$}.

A cohomology $n$--sphere over $\Z_p$ is a 
CW--complex
having the same cohomology groups with 
$\Z_p$--coefficients as the $n$--di\-men\-sio\-nal sphere $S^n$.
\begin{prop}[Volovikov's Lemma \cite{vol96}] \label{lem-vol}Set
 $G=(\Z_p)^r$, and let $X$ and $Y$ be fixed point free $G$--spaces 
such that $Y$ is a  
finite--dimensional cohomology $n$--sphere over $\Z_p$ and 
$\tilde{H}^i(X,\Z_p)=0$ for all $i\leq n$. Then there is no 
$G$--map \mbox{from $X$ to $Y$}. 
\end{prop}
Volovikov \cite{vol96} derives from this lemma a proof of the Topological
Tverberg Theorem in the prime power case. The proof of
Proposition~\ref{lem-vol} uses deeper results from bundle cohomology.
\end{section} 
\begin{section}{The extension of the lower bound}\label{sec-tow}
The next two lemmas enable us to replace the index argument used in
\cite[Section 6.6]{mat03} by Volovikov's 
Lemma. From now on let $q=p^r$ be a
prime power and $G:=(\Z_p)^r\subset S_q$ be as above.
\begin{lem}\label{lem-fixpointfree}Let $X^{*q}_{\Delta}$ be the
  $q$--fold $q$--wise deleted join for some space $X$ equipped with
the $G$--action defined as above. Then $X^{*q}_{\Delta}$ is a fixed point
free $G$--space.
\end{lem}
Note that the $G$--action on $X^{*q}_{\Delta}$ is in general not free.
\bprf Let $x=t_1x_1\oplus t_2x_2\oplus\cdots\oplus t_qx_q\in
X^{*q}_{\Delta}$, then by definition there are indices $i$ 
and $j$ such that $t_i\not=t_j$ or $x_i\not= x_j$. The indices $i$ and $j$
correspond to elements $a$ resp.~$b$ of $(\Z_p)^r$. Setting $g=b-a$,
we get $x\not=g\, x$ hence $|O_x|>1$.
\eprf
\begin{lem}\label{lem-deletedjoin}Let $q\geq 2$ and $d$ be integers. Then
we have $(\R^{d})^{*q}_{\Delta}\simeq S^{(d+1)(q-1)-1}$.
\end{lem}
\bprf Using the geometric version of the join we get an embedding\linebreak
\mbox{$(\R^d)^{*q}_{\Delta}\subset\R^{q(d+1)-1}$}. More precisely, we can identify
it with the subset\linebreak \mbox{$\{(x_1,t_1,x_2,t_2,\ldots,x_q,t_q)\,|\,
x_i\in\R^d, t_i\geq 0, \sum_1^q t_i=1\}$}. 
The diagonal is now a $d$--di\-men\-sio\-nal affine
subspace $A$, its orthogonal complement $A^\perp$ has dimension $(d+1)(q-1)$. The
restriction of the orthogonal projection $p_{A^\perp}$ onto the
complement maps $(\R^d)^{*q}_{\Delta}$ to
$\R^{(d+1)(q-1)}\setminus\{\mbox{pt}\}$. 
This map is a homotopy equivalence. 
\eprf
In the prime case, the following proof reduces to the 
Vu\'ci\'c--\v{Z}ivaljevi\'c proof, in the
version of Matou\v{s}ek \cite[Section 6.6]{mat03}.
\bprf(of {\bf  Theorem~\ref{thm-main}})
Let $\Ka$ be the simplicial complex $(\sigma^N)^{*q}_{\Delta(2)}$. 
The vertex set of $\Ka$ is $[N+1]\times [q]$.
A maximal simplex of $\Ka$ is of the form $F_1\uplus
F_2\uplus\cdots\uplus F_q$, where the $F_i$ are pairwise 
disjoint subsets of the vertex set 
$[N+1]$ of $\sigma^N$ and $\bigcup_1^q F_i=[N+1]$.
In other words, there is a one--to--one correspondence between 
the maximal simplices $\Ka$ and the ordered partitions 
$(F_1,F_2,\ldots,F_q)$ of
the vertex set $[N+1]$. Another way of looking at $\Ka$: The set of 
all maximal simplices can be identified
with the complete $(N+1)$--partite hypergraph on the vertex set $[N+1]\times [q]$.
For example, a maximal simplex in the case $d=2$ and $q=4$ encoding a
Tverberg partition for $N+1=10$ points in $\R^2$:

\begin{figure}[ht]
  \centering
  \includegraphics{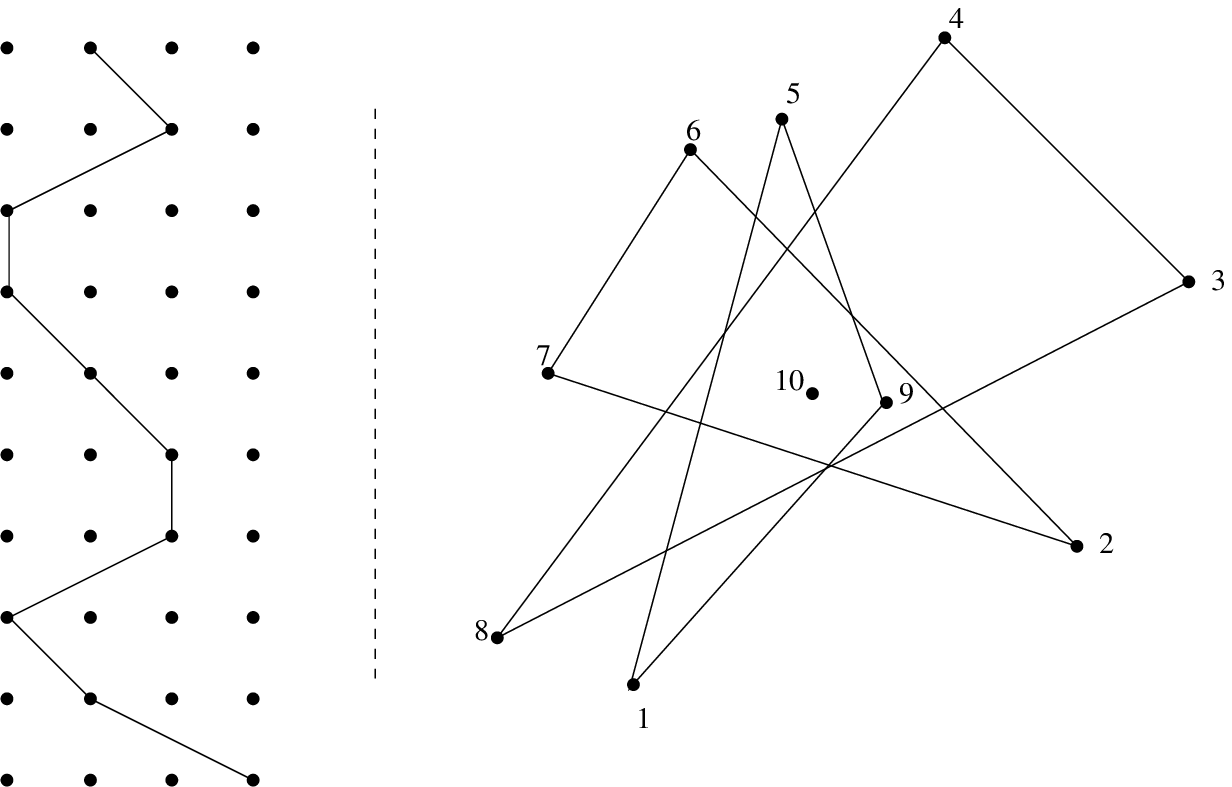}
\end{figure}

The induced $G$--action permutes the $q$ columns of vertices. 
We call a maximal face {\it good} if it encodes a
Tverberg partition of the \mbox{map $f$}. Let $f^{*q}:\|\Ka\|\ra (\R^d)^{*q}$ be
the $q$--fold join of $f$ restricted to $\|\Ka\|$, then it is a $G$--map. 
A maximal simplex $S$ of $\Ka$ is good if its image
$f^{*q}(\| S\|)$ intersects the diagonal of $(\R^d)^{*q}$. 
Proving a lower bound for the number of good simplices in 
$\Ka$ gives then a lower bound for the number of Tverberg partitions of $f$.
If there are at least $M$ good simplices we have a least
$M/q!$ unordered Tverberg partitions.

In the next paragraph, we define a family $\Le$ of subcomplexes $\La\subset \Ka$
having the properties: (i) $\La$ is closed under
the $G$--action, and (ii) conn$(\La)\geq N-1$. Then $\La$ is 
again a fixed point free $G$--space by (i) and Lemma~\ref{lem-fixpointfree}.
The reduced cohomology groups of $\La$ vanish in
dimensions $0$ to $N-1$ due to (ii). 
Now with Lemma~\ref{lem-deletedjoin} we get as a
direct corollary of Volovikov's Lemma that $\La$ contains one good maximal
\mbox{simplex $S$}; in fact, the entire orbit of $S$ is good and we get $q$
good simplices in $\La$. Suppose $Q$ is the number of ${\sf
  L}\in\Le$ containing any given maximal simplex of $\Ka$, then we
obtain the lower bound
\begin{eqnarray}\label{bound}
M\geq q\cdot |\Le |/Q.
\end{eqnarray}

We define the family $\Le$ and distinguish two cases: (i) $N$ even,
that is, $p$ or $d$ is odd, and (ii) $N$ odd, that is,
$p=2$ and even $d$. First we 
divide the $N+1$ rows into pairs
such that we get $\frac{N}{2}$ pairs and one remaining row in the
first case, and $\frac{N+1}{2}$ pairs in the second. Now we focus on
the two rows of one pair; the simplices of $\Ka$ living on these
two rows form bipartite graphs $K_{q,q}$. Suppose that we have chosen 
a connected $G$--invariant subgraph $C_i$ of $K_{q,q}$, 
$i\in[\frac{N}{2}]$ resp.~$i\in [\frac{N+1}{2}]$, for every pair. 
The maximal simplices of $\La$ to a given 
choice of row pairing and of the $C_i$, $i\in[\frac{N}{2}]$ resp.~$i\in
[\frac{N+1}{2}]$, are the maximal simplices of $\Ka$ that contain an 
edge of each $C_i$. $\La$ is $G$--invariant by construction. 
Topologically, we get in the first case
\[ \La=C^{*(N/2)}*D_q,\]
and in the second
 \[ \La=C^{*((N+1)/2)}.\]
Here $D_q$ is the discrete space on $q$ elements; in both cases one 
has $\mbox{conn}(\La)\geq N-1$ using inequality~(\ref{conn}).

Now we explain how to get 
\[q(q-p^0)(q-p^1)(q-p^2)\cdots(q-p^{r-1})/(r+1)! \]
distinct $G$--invariant, 
connected subgraphs $C$ by choosing $r+1$ edges of $K_{q,q}$. 
For $q$ prime, this process coincides with the construction described
in \cite[Section 6.6]{mat03}.
To obtain a $G$--invariant subgraph choose edges and take their orbits,
see Figure~\ref{fig-orbit} for orbits in the case $q=3^2$. 
The vertices are elements of $(\Z_p)^r$ 
having order $p$ as group elements. To make sure that we count an 
orbit without multiplicities choose its representative edge as 
the edge that is incident to the upper left vertex $O:=(0,0,\ldots,0)$.

\begin{figure}[ht]

    \centering
  \includegraphics{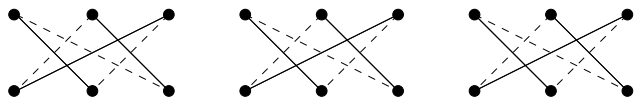}
  \caption{$G$--orbits of the edges $((0,0),(0,1))$ and 
    $((0,0),(0,2))$.}
  \label{fig-orbit}
\end{figure}

To prove the connectivity of the graph $C$ we show that the
component $K_O$ of the vertex $O$ is the whole graph $C$.
Choosing $r+1$ representative edges consecutively such that
in each step a new component is connected to the component
$K_O$ leads to a connected subgraph. 

More precisely, we will show inductively that 
after $1\leq k\leq r+1$ steps: (i) there are
$2p^{k-1}$ vertices in each component,
$p^{k-1}$ in each shore, and (ii) in total there
are $p^{r-(k-1)}$ components. For $k=1$, the orbit of 
an edge consists of $p^r$ vertex--disjoint edges, see Figure~\ref{fig-orbit}. 
For $k=2$, the graph of two orbits is equal to the disjoint union
of $p^{r-1}$ cycles of length $2p$, see Figure~\ref{fig-orbit}. 
Assume that for
$1\leq k\leq r$ edges the statement is true. 
Let the $(k+1)$--st edge be an edge connecting $K_O$
with one of the other remaining $p^{r-(k-1)}-1$ components,
there are $r-p^{k-1}$ many representative edges to do so. 
The graph of the $(k+1)$--st orbit and any of
the $k$ first orbits is again a union of cycles of length $2p$,
hence each $p$ components of the graph of the first
$k$ orbits get connected. 
Therefore the number of components decreases by a factor $p$,
and the number of vertices increases by the factor $p$
in each shore.


As the order in the $r+1$ steps of our construction does not play any
role this process leads to the desired number of graphs $C$.
Every given edge determines an orbit, hence there are   
\[(q-p^0)(q-p^1)(q-p^2)\cdots(q-p^{r-1})/r!\]
connected, $G$--invariant graphs $C$ containing this edge.

Finally, let $\pi$ be the number of possibilities to do
the row pairing in \mbox{case (i) or (ii)} ($\pi$ cancels
out in the end). Then in case (i) we get:
\begin{eqnarray*}
& |\Le |= \pi\cdot\left( q\cdot\prod\nolimits_{i=0}^{r-1}(q-p^i)/(r+1)!
\right) ^{N/2},\\
& Q=\pi\cdot\left(\prod\nolimits_{i=0}^{r-1}(q-p^i)/r!\right)^{N/2},
\end{eqnarray*}
and in case (ii):
\begin{eqnarray*}
& |\Le |= \pi\cdot\left(q\cdot\prod\nolimits_{i=0}^{r-1}(q-p^i)
/(r+1)!\right)^{(N+1)/2},\\
& Q=\pi\cdot\left(\prod\nolimits_{i=0}^{r-1}(q-p^i)
  /r!\right)^{(N+1)/2}.
\end{eqnarray*}
Plugging these numbers into inequality~(\ref{bound}) completes the proof.
\eprf
\end{section}
\begin{section}{On the number of splitting necklaces}\label{sec-necklace}
It is known that the methods introduced for the 
Topological Tverberg Theorem can also be applied 
to the splitting problem for necklaces for many thieves,
see \cite[Section 6.4]{mat03}. We will extend the lower bound
of \cite{vz93} to the prime power case.
A {\it necklace} is modeled in the 
following way: Given $d$ continuous probability measures
on $[0,1]$ and $q\geq 2$ thieves. A {\it fair splitting}
of the necklace consists of a partition of $[0,1]$
into a number $n$ of subintervals $I_1,I_2,\ldots,
I_n$ and a partition of $[n]$ into $q$ subsets 
$T_1,T_2,\ldots ,T_q$ such that every thief 
has an equal amount of all $d$ materials: 
\[ \sum\nolimits_{j\in T_k}\mu_i(I_j)=\frac{1}{q}\mbox{
, for all $1\leq i\leq d$ and $1\leq k\leq q$}.\]
Noga Alon proved in 1987 that in general $d(q-1)$ is 
the smallest number of cuts for $q$ thieves. A
necklace is called {\it generic} if there is 
no fair splitting with less than $d(q-1)$ cuts.
The following result extends the lower bound of \cite{vz93}
for the number of fair splittings to the prime power
case.
\begin{thm}\label{thm-necklace}Let $q=p^r$ be a prime
power. For generic necklaces made out of $d$ continuously 
distributed materials the number of fair splittings with
$d(q-1)$ cuts for $q$ thieves is at least:
\[q\cdot\left(\frac{q}{r+1}\right)^{\lceil\frac{d(q-1)}{2}\rceil}. \]
\end{thm}
In the proof we will again face deleted joins, but also
the deleted product $(\R^d)^q_\Delta$ that is  
the q--fold cartesian product of $\R^d$ without its diagonal. It is
well--known that $(\R^d)^q_\Delta\simeq S^{d(q-1)-1}$, 
see e.~g.~\cite[Section 6.3]{mat03}.
\bprf(sketch) In the proof of Theorem 6.4.1 of \cite{mat03}
there is a one--to--one correspondence between the set of
splittings  of a generic necklace for $q$ thieves and the simplicial
complex $\Ka=(\sigma^{d(q-1)+1})^{*q}_{\Delta(2)}$. The map 
$f:\|\Ka \|\ra (\R^d)^q,\,z\mapsto f(z)_{i,k}:=\sum\nolimits_{j\in T_k}
\mu_i(I_j)$ expressing the gains of the thieves is a $G$--map.
If there is no fair splitting, $f$ would miss the diagonal
of $(\R^d)^q$. Now let $\Le$ be a family of subcomplexes $\La$
satisfying: (i) $\La$ is closed under the $G$--action, and (ii) 
conn$(\La)\geq d(q-1)-1$. Again with Volovikov's Lemma
every $\La$ contains at least one fair splitting, but as above
the whole orbit of size $q$ is good. In conclusion, the whole 
construction for $\Le$ and the counting as in the proof of
Theorem~\ref{thm-main} can be carried over.
\eprf
\end{section}



\begin{thebibliography}{99999999}


\bibitem[BSS81]{bss81}I.~B\'ar\'any, S.~B.~Shlosman and A.~Sz\"ucs, {\it On a topological generalization of a theorem of Tverberg}, J.~London Math.~Soc., \mbox{II.~Ser.~23 (1981), 158--164}


\bibitem[Mat03]{mat03}J.~Matou\v{s}ek, {\it Using the Borsuk--Ulam Theorem}, Springer (2003)

\bibitem[\"Oz87]{oez87}M.~\"Ozaydin, {\it Equivariant maps for the symmetric group}, Preprint, University of Wisconsin--Madison, 1987, 17 pages



\bibitem[Sch04]{sch04}T.~Sch\"oneborn, {\it On the Topological
    Tverberg Theorem}, diploma thesis, TU Berlin (2004)

\bibitem[Tve66]{tve66}H.~Tverberg, {\it A generalization of Radon's Theorem}, J.~London Math.~Soc., 41 (1966), 123--128

\bibitem[Vol96]{vol96}A.~Yu.~Volovikov, {\it On a topological generalization of the Tverberg Theorem}, Math.~Notes, 3 (1996), 324--326

\bibitem[V\v{Z}93]{vz93}A.~ Vu\'ci\'c and R.~\v{Z}ivaljevi\'c, {\it Notes on a conjecture of Sierksma},
Discr.~Comput.~Geom., 9 (1993), 339--349



\end{thebibliography}
\end{document}